\documentstyle[amsmath,amssymb,12pt]{article}

\newcommand{\rom}[1]{{\rm #1}}

\allowdisplaybreaks

\begin{document}

\newtheorem{definition}{Definition}
\newtheorem{remark}{Remark}
\newtheorem{proposition}{Proposition}
\newtheorem{theorem}{Theorem}
\newtheorem{corollary}{Corollary}
\newtheorem{lemma}{Lemma}

\newcommand{\indlim}{\operatornamewithlimits{ind\,lim}}
\newcommand{\Ffin}{{\cal F}_{\mathrm fin}}
\newcommand{\Fext}{{\cal F}_{\mathrm ext}}
\newcommand{\D}{{\cal D}}
\newcommand{\N}{{\Bbb N}}
\newcommand{\C}{{\Bbb C}}
\newcommand{\Z}{{\Bbb Z}}
\newcommand{\R}{{\Bbb R}}
\newcommand{\Rp}{{\R_+}}
\newcommand{\eps}{\varepsilon}
\newcommand{\supp}{\operatorname{supp}}
\newcommand{\la}{\langle}
\newcommand{\ra}{\rangle}
\newcommand{\const}{\operatorname{const}}
\renewcommand{\emptyset}{\varnothing}
\newcommand{\di}{\partial}
\newcommand{\hotimes}{\hat\otimes}

\newcommand{\pii}{\pi_{\nu\otimes\sigma}}
\newcommand{\RR}{{\cal R}}
\newcommand{\RX}{{\RR\times X}}
\newcommand{\ZZ}{\Z_{+,\,0}^\infty}

\begin{center}{\Large \bf The Jacobi field of a L\'evy process}\end{center}

{\large Yuri M. Berezansky}

\noindent{\sl Institute of  Mathamatics NAN Ukraine,
Tereshchenkivska 3, 01601 Kiev, Ukraine}\\[2mm] {\rm E-mail:
berezan@mathber.carrier.kiev.ua} \vspace{5mm}

{\large Eugene Lytvynov}

\noindent{\sl Institut f\"{u}r Angewandte Mathematik,
Universit\"{a}t Bonn, Wegelerstr.~6, D-53115 Bonn,\\ Germany;
BiBoS, Univ.\ Bielefeld, Germany\\[2mm]
 {\rm E-mail: lytvynov@wiener.iam.uni-bonn.de}}\vspace{5mm}

{\large Dmytro A. Mierzejewski}

\noindent{\sl Zhytomyr State Pedagogical University, Velyka
Berdychivska Street 40, 10008 Zhytomyr,  Ukraine}\\[2mm] {\rm
E-mail: mierz@zspu.edu.ua}

\begin{abstract}

\noindent We derive an explicit formula for the Jacobi field that
is acting in an extended Fock space and corresponds to an
($\R$-valued) L\'evy process on a Riemannian manifold. The support
of the measure of jumps in the L\'evy--Khintchine representation
for the L\'evy process is supposed to have an infinite number of
points. We characterize the
 gamma,
Pascal, and Meixner processes as the only L\'evy processes whose
Jacobi field leaves the set of finite continuous elements of the
extended Fock space invariant.

\end{abstract}

\noindent 2000 {\it AMS Mathematics Subject Classification}.
Primary: 60G51, 60G57, 47B36. Secondary: 60H40.\vspace{3mm}

The aim of this notice is to derive an explicit formula for the
 Jacobi field \cite{BeLi,Ly, bere,new,Chap} that is
acting in an extended Fock space \cite{ silva,KL, beme1,Ly3,Ly4}
and corresponds to an ($\R$-valued) L\'evy process on a Riemannian
manifold. The support of the measure of jumps in the
L\'evy--Khintchine representation for the process is supposed to
have  an infinite number of points. The proof of this formula will
be based on a result of \cite{Ly4}, see also \cite{NS,S}. We will
characterize the  gamma, Pascal, and Meixner processes as the only
L\'evy processes whose Jacobi field leaves the set of finite
continuous elements of the extended Fock space invariant.

So, let $X$ be a complete, connected, oriented $C^\infty$
(non-compact) Riemannian manifold and let ${\cal B}(X)$ be the
Borel $\sigma$-algebra on $X$. Let $\sigma$ be a Radon measure on
$(X,{\cal B}(X))$ that is non-atomic and non-degenerate (i.e.,
$\sigma(O)>0$ for any open set $O\subset X$). We denote by $\D$
the space $C_0^\infty(X)$ of all infinitely differentiable,
real-valued functions on $X$ with compact support. It is known
that $\D$ can be endowed with a topology such that the natural
embedding of $\D$ into the real $L^2$-space $L^2(X;\sigma)$ is
dense, continuous, and nuclear. Thus, we can consider the standard
nuclear triple $ \D'\supset L^2(X;\sigma)\supset\D$,  where $\D'$
is the dual space of $\D$ with respect to the zero space
$L^2(X;\sigma)$. The dual pairing between $\omega\in\D'$ and
$f\in\D$  will be denoted by $\la\omega,f\ra$.  By ${\cal C}(\D')$
we will denote the cylinder $\sigma$-algebra on $\D'$.

Let  ${\cal R}{:=}\R\setminus\{0\}$. We endow $\RR$ with the
relative topology of $\R$ and let ${\cal B}(\RR)$ be the Borel
$\sigma$-algebra on $\RR$. Let $\nu$ be a Radon measure on
$(\RR,{\cal B}(\RR))$ whose  support contains an infinite number
of points. Let $\tilde\nu(ds){:=}s^2\,\nu(ds)$. We suppose that
$\tilde \nu$ is a finite measure on $(\RR,{\cal B}(\RR))$, and
furthermore,   there exists $\varepsilon>0$
 such that
\begin{equation}\label{4335ew4} \int_{\RR} \exp\big(\varepsilon
|s|\big)\,\tilde\nu(ds)<\infty.\end{equation}  Therefore, the
measure $\tilde\nu$ has all moments finite, and  the set of all
polynomials is  dense in $L^2(\RR;\tilde\nu)$.

 We now
define a centered L\'evy process on $X$  as a generalized process
on ${\cal D}'$ whose law is the probability measure
$\rho_{\nu,\,\sigma} $ on $({\cal D}',{\cal C}({\cal D}'))$ given
by its Fourier transform
 \begin{equation}\label{rew4w}
\int_{{\cal D}'} e^{i\la \omega,\varphi\ra}\,
\rho_{\nu,\,\sigma}(d\omega)=\exp\bigg[\int_\RX(e^{is\varphi(x)}-1-is\varphi(x))\,\nu(ds)\,\sigma(dx)\bigg],\qquad
\varphi\in {\cal D}\end{equation} (compare with \cite{GV, Tsil}).
The existence of $\rho_{\nu,\,\sigma}$ follows from the
Bochner--Minlos theorem. Formula \eqref{rew4w} is the
L\'evy--Khintchine representation for the L\'evy process.

We will now construct a decomposition of  the $L^2$-space
$L^2(\D';\rho_{\nu,\,\sigma})$ following the idea of
orthogonalization of continuous polynomials with respect to a
probability measure  that is defined on a co-nuclear space, cf.\
\cite[Sect.~11]{Sko}.

We denote by ${\cal P}({\cal D}')$ the set of continuous
polynomials on ${\cal D}'$, i.e., functions on ${\cal D}'$ of the
form $F(\omega)=\sum_{i=0}^n\la\omega^{\otimes i},f_{i}\ra$, $
\omega^{\otimes 0}{:=}1,\ f_{i}\in{\cal D}^{\hotimes i}$,
$i=0,\dots,n$, $n\in\Z_+$.  Here, $\hat\otimes$ stands for
symmetric tensor product.  The greatest number $i$ for which
$f^{(i)}\ne0$ is called the power of a polynomial. We denote by
${\cal P}_n({\cal D}')$ the set of continuous polynomials of power
$\le n$.

By \eqref{4335ew4}, \eqref{rew4w}, and \cite[Sect.~11]{Sko},
${\cal P}({\cal D}')$ is a dense subset of $ L^2({\cal
D}';\rho_{\nu\otimes\sigma})$.
Let ${\cal P}^\sim _n({\cal D}')$ denote the closure of ${\cal
P}_n({\cal D}')$ in $L^2({\cal D}';\rho_{\nu,\,\sigma})$, let
${\bf P}_n({\cal D}')$, $n\in\N$, denote the orthogonal difference
${\cal P }^\sim_n({\cal D}')\ominus{\cal P}^\sim_{n-1}({\cal
D}')$, and let ${\bf P}_0({\cal D}'){:=}{\cal P }^\sim_0({\cal
D}')$. Then, we evidently have: \begin{equation}\label{zgzguzug}
L^2({\cal D}';\rho_{\nu,\,\sigma})=\bigoplus_{n=0}^\infty{\bf
P}_n({\cal D}'). \end{equation}

The set of all projections ${:}\la\cdot^{\otimes n},f_n\ra{:}$ of
continuous monomials $\la \cdot^{\otimes n},f_n\ra$,
$f_n\in\D^{\hat\otimes n}$, onto ${\bf P}_n({\cal D }')$ is dense
in ${\bf P}_n({\cal D}')$. For each $n\in\N$, we define a Hilbert
space ${\frak F}_n$
 as the closure of the set ${\cal D}^{\hotimes n}$ in the
 norm generated by the scalar product \begin{equation}\label{gffdz}
 (
 f_n,g_n)_{{\frak F}_n}{:=}\frac1{n!}\, \int_{{\cal D}'}{:}\la \omega^{\otimes n},f_n\ra{:}
 \, {:}\la\omega^{\otimes n},g_n\ra{:}\,\rho_{\nu,\,\sigma}(d\omega),\qquad f_n,g_n\in{\cal D}^{\hotimes n}.
 \end{equation} Denote  \begin{equation}\label{qztztzt}{\frak
F}{:=}\bigoplus_{n=0}^\infty{\frak F}_n\,n!,\end{equation} where
${\frak F}_0{:=}\R$. By \eqref{zgzguzug}--\eqref{qztztzt}, we get
 the unitary operator ${\cal U}:{\frak F}\to
L^2(\D';\rho_{\nu,\,\sigma})$ that is defined through $ {\cal
U}f_n{:=}{:}\la\cdot^{\otimes n},f_n\ra{:}$, $f_n\in{\cal
D}^{\hotimes n}$, $n\in\Z_+$, and then extended by linearity and
continuity to the whole space ${\frak F}$\rom.

We will now write down an explicit formula for the scalar product
$(\cdot,\cdot)_{{\frak F}_n}$. In the case of the Gamma process,
this formula is due to  \cite{KL} and  \cite{silva}, in the case
of the Pascal and Meixner process due to \cite{Ly3} and
\cite{bere3}, and in the case of a general L\'evy process due to
\cite{Ly4}.

We denote by $\ZZ $ the set of all sequences $\alpha$ of the form
$\alpha=(\alpha_1,\alpha_2,\dots,\alpha_n,0,0,\dots)$,
$\alpha_i\in\Z_+$, $n\in\N$. Let $|\alpha|{:=}\sum_{i=1}^\infty
\alpha_i$. For each $\alpha\in\ZZ$, $1\alpha_1+2\alpha_2+\dots=n$,
$n\in\N$, and for any function $f_n:X^n\to\R$ we define a function
$D_\alpha f_n:X^{|\alpha|}\to\R$ by setting \begin{align}(D_\alpha
f_n)(x_1,\dots,x_{|\alpha|}){:=}& f(x_1,\dots,x_{\alpha_1},
\underbrace{x_{\alpha_1+1},x_{\alpha_1+1}}_{\text{2 times }},
\underbrace{x_{\alpha_1+2},x_{\alpha_1+2}}_{\text{2 times
}},\dots,
\underbrace{x_{\alpha_1+\alpha_2},x_{\alpha_1+\alpha_2}}_{\text{2
times }},\notag\\ &\quad
\underbrace{x_{\alpha_1+\alpha_2+1},x_{\alpha_1+\alpha_2+1},,x_{\alpha_1+\alpha_2+1}}_{\text{3
times }},\dots).\label{123}\end{align} Let $(
P_n(\cdot))_{n=0}^\infty$ be the system of polynomials  with
leading coefficient 1 that are orthogonal with respect to the
measure $\tilde\nu(ds)$ on $\RR$. We have,  for any
$f_n,g_n\in{\cal D}^{\hotimes n}$\rom, $n\in\N$\rom,
\begin{gather}(f_n,g_n)_{{\frak F}_n} =\sum_{\alpha\in\ZZ:\,
1\alpha_1+2\alpha_2+\dots=n}K_\alpha \int_{X^{|\alpha|}}(D_\alpha
f_n)(x_1,\dots,x_{|\alpha|})\notag
\\ \times (D_\alpha g_n)(x_1,\dots,x_{|\alpha|}) \,\sigma^{\otimes
|\alpha|}(dx_1,\dots,dx_{|\alpha|}),\label{atfztftfz}\end{gather}
where \begin{equation}\label{tfztftfz} K_\alpha=
\frac{n!}{\alpha_1!\,\alpha_2!\dotsm}\,\prod_{k\ge1}\bigg(\frac{\|P_{k-1}\|_{L^2(\RR;\tilde\nu)}}{k!}\bigg)^{2\alpha_k}.
\end{equation}

By using \eqref{atfztftfz}, \eqref{tfztftfz}, one derives the
following representation of ${\frak F}_n$ (see
\cite[formula~(5.19)]{Ly4}):
\begin{equation}\label{jhjhhj}{\frak F
}_n{:=}\bigoplus_{\alpha\in\ZZ:\,
1\alpha_1+2\alpha_2+\dots=n}{\frak F}_{n,\,\alpha},\qquad {\frak
F}_{n,\,\alpha}{:=}
L^2_\alpha(X^{|\alpha|};\sigma^{\otimes|\alpha|})\,K_\alpha.\end{equation}
Here, $$L^2_{\alpha}(X^{|\alpha|};\sigma^{\otimes |\alpha|})=
L^2(X;\sigma)^{\hotimes \alpha_1}\otimes L^2(X;\sigma)^{\hotimes
\alpha_2}\otimes\dotsm\,,$$ and for each $f_n\in{\cal
D}^{\hat\otimes n}\subset{\frak F}_n$, the ${\frak
F}_{n\,,\alpha}$-coordinate of $f_n$ is equal to $D_\alpha f_n$.
Thus, we can extend $D_\alpha$ by continuity to the orthogonal
projection of ${\frak F}_n$ onto ${\frak F}_{n\,,\alpha}$. In what
follows, we will also denote by $S_\alpha$ the orthogonal
projection of $L^2(X^{|\alpha|};\sigma^{\otimes|\alpha|})$ onto
$L_\alpha^2(X^{|\alpha|};\sigma^{\otimes|\alpha|})$.  Taking
\eqref{jhjhhj} into account, we will call $\frak F$ an extended
Fock space (compare with \cite{KL, beme1}).

For an arbitrary  $\varphi\in\D$, we consider in the space
$L^2(\D';\rho_{\nu,\,\sigma})$ the operator $M(\varphi)$ of
multiplication by the function $\la\cdot,\varphi\ra$, and let
$J(\varphi){:=}{\cal U}M(\varphi){\cal U}^{-1}$. We denote by
$\Ffin(\D)$ the set of all vectors of the form
$(f_0,f_1,\dots,f_n,0,0,\dots)$, $f_i\in\D^{\hotimes i}$,
$i=0,\dots,n$, $n\in\Z_+$. Evidently, $\Ffin(\D)$ is a dense
subset of $\frak F$.

\begin{theorem}\label{zgztpopo}
For any $\varphi\in\D$\rom, we have\rom:
\begin{equation}\label{7t667}
\Ffin(\D)\subset\operatorname{Dom}(J(\varphi)),\qquad
J(\varphi)\restriction\Ffin(\D)=J^+(\varphi)+J^0(\varphi)+J^-(\varphi),\end{equation}
 the linear operators  $J^+(\varphi)$\rom,
$J^0(\varphi)$\rom, $J^-(\varphi)$ being defined as follows\rom:
for any $f_n\in\D^{\hotimes n}$\rom, $n\in\Z_+$\rom,
\begin{equation}\label{hbgg} J^+(\varphi) f_n=\varphi\hotimes
f_n,\end{equation} $J^0(\varphi)f_n\in {\frak F}_n$ and each
${\frak F}_{n,\,\alpha}$-coordinate of $J^0(\varphi)f_n$ is equal
to
\begin{equation}\label{gfg}
\big(J^0(\varphi)f_n\big)_\alpha(x_1,\dots,x_{|\alpha|})=\sum_{k=1}^\infty
\alpha_k a_{k-1}
S_\alpha\big(\varphi(x_{\alpha_1+\dots+\alpha_k})(D_\alpha
f_n)(x_1,\dots,x_{|\alpha|})\big)\qquad
\text{$\sigma^{\otimes|\alpha|}$-a.e.,} \end{equation}
$J^-(\varphi)f_n=0$ if $n=0$\rom, $J^-(\varphi)f_n\in {\frak
F}_{n-1}$ if $n\in\N$\rom, and each ${\frak
F}_{n-1,\,\alpha}$-coordinate of $J^-(\varphi)f_n$ is equal to
\begin{multline}
\big(J^-(\varphi)f_n\big)_\alpha(x_1,\dots,x_{|\alpha|})=n\tilde\nu({\cal
R }) S_\alpha\bigg(\int_X
\varphi(x)(D_{\alpha+1_1}f_n)(x,x_1,\dots,x_{|\alpha|})\,\sigma(dx)\bigg)\\\text{}+\sum_{k\ge
2 }\frac nk\, \alpha_{k-1}b_{k-1} S_\alpha \big(
\varphi(x_{\alpha_1+\dots+\alpha_k})(D_{\alpha-1_{k-1}+1_k}f_n)(x_1,\dots,x_{|\alpha|})\big)\qquad
\text{$\sigma^{\otimes|\alpha|}$-a.e.} \label{ggztr}\end{multline}
In formulas \eqref{gfg} and \eqref{ggztr}\rom, we denoted $
\alpha\pm
1_n{:=}(\alpha_1,\dots,\alpha_{n-1},\alpha_n\pm1,\alpha_{n+1},\dots)$
for $\alpha\in \ZZ$ and $n\in\N$\rom;  the real numbers $a_n$ and
positive numbers $b_n$ are given through the recurrence relation
\begin{equation}\label{12345} s P_n(s)=P_{n+1}(s)+a_n
P_n(s)+b_n P_{n-1}(s),\qquad n\in\Z_+,\  P_{-1}(s){:=}0.
\end{equation} Finally\rom, $J(\varphi)$ is essentially
self-adjoint on $\Ffin(\D)$\rom.

\end{theorem}

By \eqref{7t667}, the operator $J(\varphi)\restriction\Ffin(\D)$
is a sum of creation, neutral, and  annihilation  operators, and
hence $J(\varphi)\restriction\Ffin(\D)$ has a Jacobi operator
structure. The family of operators $(J(\varphi))_{\varphi\in \D}$
is called the Jacobi field corresponding to the L\'evy processes
with law $\rho_{\nu,\,\sigma}$.

The proof of \eqref{7t667}--\eqref{ggztr} follows from
\cite[Theorem~5.1, Corollaries~4.2 and~5.1]{Ly4}. The essential
self-adjointness of $J(\varphi)$ on $\Ffin(\D)$ follows from
\eqref{7t667}, \eqref{hbgg} and \cite[Theorem~4.1]{new} whose
proof admits a direct generalization to the case of the extended
Fock space $\frak F$.

We notice that the operator $J^+(\varphi)$ leaves the set
$\Ffin(\D)$ invariant, while the operators $J^0(\varphi)$ and
$J^-(\varphi)$, in general do not.

\begin{corollary}\label{bfgv}
Suppose that\rom, for each $\varphi\in\D$\rom,
$J^0(\varphi)\Ffin(\D)\subset \Ffin(\D)$ and
$J^-(\varphi)\Ffin(\D)\subset \Ffin(\D)$\rom, so that
$J(\varphi)\Ffin(\D)\subset \Ffin(\D)$\rom. Then\rom, $\tilde\nu$
is a finite measure on $\cal R$ such that $a_n=\lambda(n+1)$ and
$b_n=\varkappa n(n-1)$\rom, $n\in\Z_+$\rom. Here\rom, $a_n$ and
$b_n $ are the coefficients from
 \eqref{12345}\rom, and $\lambda\in\R$ and $\varkappa>0$ are
 arbitrarily chosen parameters\rom. Furthermore\rom, we have in
 this case, for each $f_n\in\D^{\hotimes n}$, $n\in\Z_+$:
\begin{align*} &\big(J^0(\varphi)f_n\big)(x_1,\dots,x_n)=\lambda
n\big(\varphi(x_1)f_n(x_1,\dots,x_n)\big)\hat{},\qquad
\text{$\sigma^{\otimes n}$-a.e.,}\\
&\big(J^-(\varphi)f_n\big)(x_1,\dots,x_{n-1})=n \tilde \nu({\cal R
}) \int_{\cal R} \varphi(x)f_n(x,x_1,\dots,x_{n-1})\,\sigma(dx)\\
&\qquad +\varkappa n(n-1)
\big(\varphi(x_1)f_n(x_1,x_1,x_2,x_3,\dots,x_{n-1})\big)\hat{},\qquad\text{$\sigma^{\otimes(n-1)}$-a.e.}
\end{align*}
Here\rom, $(\cdot)\hat{}$ denotes symmetrization of a
function\rom. The choice $|\lambda|=2$ corresponds to a gamma
process\rom, $|\lambda|>2$ corresponds to a Pascal process\rom,
and $|\lambda|<2$ corresponds to a Meixner process\rom.

\end{corollary}

 Corollary~\ref{bfgv} is  derived from
Theorem~\ref{zgztpopo} and \cite[Corollary~5.1]{Ly4}, by using the
idea that that the off-diagonal values of a continuous function of
several variables uniquely determine the on-diagonal values of
this function.

The Jacobi fields of the gamma, Pascal, and Meixner processes (in
the case $\tilde\nu({\cal R})=1$ and $\varkappa=1$) were studied
in \cite{Ly3,Ly4}, see also \cite{bere3,KL,NS,S}. \vspace{2mm}

 \noindent {\bf Acknowledgements.} The first named author was partially
 supported by INTAS, Project 00-257
 and the DFG, Project 436 UKR 113/61.   The second
author acknowledges the financial support of the SFB 611, Bonn
University, and  the DFG Research Project 436 UKR 113/43.


\begin{thebibliography}{99}


\bibitem{BeLi} Yu.~M. Berezansky, V. O. Livinsky and E. W. Lytvynov,  A
  generalization of Gaussian white noise analysis, {\it Meth.\ Func.\ Anal.\ and
  Topol.}\ {\bf 1} (1995), no.~1, 28--55.

\bibitem{Ly} E.  Lytvynov, Multiple Wiener integrals and non-Gaussian white
  noises: a Jacobi field approach, {\it Meth.\ Func.\ Anal.\ and Topol.}\ {\bf
  1} (1995), no.~1, 61--85.

\bibitem{bere} Yu.\ M. Berezansky,  Commutative Jacobi fields in Fock space,
{\it Integral Equations Operator Theory\/} {\bf 30} (1998),
163--190.

\bibitem{new} Yu.\ M. Berezansky, On the theory of commutative
Jacobi fields, {\it Methods Funct.\ Anal.\ Topology\/} {\bf 4}
(1998), no.~1, 1--31.

\bibitem{Chap} Yu.\ A. Chapovsky, On the inverse spectral problem
for a commutative field of operator-valued Jacobi matrices,  {\it
Methods Funct.\ Anal.\ Topology} {\bf 8} (2002), no.~1, 14--22.


\bibitem{silva} Yu.\ G. Kondratiev, J. L. Silva, L. Streit and G. F. Us,
  Analysis on Poisson and Gamma spaces, {\it Infin.\ Dimen.\ Anal.\
   Quant.\  Probab.\ Rel.\ Top.}\ {\bf 1} (1998), 91--117.



\bibitem{KL} Y. Kondrtatiev and E. Lytvynov,  Operators of gamma
white noise calculus, {\it Infin.\ Dimen.\ Anal.\ Quant. Prob.\
Rel.\ Top.}\ {\bf 3} (2000), 303--335.


\bibitem{beme1} Yu.~M.~Berezansky and D. A. Mierzejewski,
The structure of the extended symmetric Fock space, {\it Methods
Funct.\ Anal.\ Topology\/} {\bf 6} (2000), no.~4, 1--13.

\bibitem{Ly3} E. Lytvynov,  Polynomials of Meixner's type in
infinite dimensions---Jacobi fields and orthogonality measures,
{\it J. Funct.\ Anal.}, to appear.

\bibitem{Ly4} E. Lytvynov,  Orthogonal decompositions for L\'evy processes with an application to the
gamma, Pascal, and Meixner processes, {\it Infin.\ Dimen.\ Anal.\
Quant. Prob.\ Rel.\ Top.}, to appear.


\bibitem{NS}  D. Nualart and W. Schoutens, Chaotic and
predictable representations for L\'evy processes,  {\it Stochastic
Process.\ Appl.}\ {\bf 90} (2000),  109--122.


\bibitem{S} W. Schoutens, ``Stochastic Processes and Orthogonal
Polynomials,'' Lecture Notes in Statist., Vol.~146,
Springer-Verlag, New York, 2000.

\bibitem{GV} I. M. Gel'fand and N. Ya.~Vielenkin, ``Generalized
Functions, Vol.~4. Applications of Harmonic Analysis,'' Academic
Press,  New York/London, 1964.

 \bibitem{Tsil} N. Tsilevich, A. Vershik, and M. Yor,   An
infinite-dimensional analogue of the Lebesgue measure and
distinguished properties of the gamma process, {\it J.\ Funct.\
Anal.}\ {\bf 185 } (2001), 274--296.

\bibitem{Sko} A. V. Skorohod, ``Integration in Hilbert Space,''
Springer-Verlag, New York, 1974.



\bibitem{bere3} Yu.~M. Berezansky,  Pascal measure on generalized functions and the corresponding
generalized Meixner polynomials,  {\it  Methods Funct.\ Anal.\
Topology} {\bf 8} (2002), no.~1, 1--13.















\end{thebibliography}
\end{document}